\documentclass[fleqn,10pt]{article}

    \usepackage{amsmath,amssymb,amsfonts,amsbsy,amsthm}
    \usepackage{latexsym,graphics,graphicx,color,sectsty,mathdots}
    
    \usepackage{subcaption}
    \usepackage{layout}
    \usepackage{bbm}
    \usepackage[colorlinks=true,linkcolor=blue]{hyperref}
    \usepackage{calc,cancel}
    \usepackage{enumitem}
    \usepackage{rotating}
    \usepackage[table]{xcolor}
    \usepackage{cleveref}
    \usepackage{appendix} 
    \usepackage{sidecap}
    \usepackage{mathtools}

    \usepackage{arydshln}
    \usepackage{makeidx} 

\theoremstyle{plain}

\theoremstyle{remark}
\newtheorem{example}{Example}

	\definecolor{bgblue}{rgb}{0.04,0.39,0.53}
	\definecolor{dblue}{rgb}{0,0.3,0.7}
	\definecolor{ddblue}{rgb}{0,0.1,0.6}
	\definecolor{ddgreen}{rgb}{0,0.25,0.05}
	\definecolor{dgreen}{rgb}{0,0.5,0.05}

\newcommand{\enma}[1]   {\ensuremath{#1}}

\newcommand{\req}[1]{(\ref{#1.eq})}

\newcommand{\beq}{\begin{equation}}
\newcommand{\eeq}{\end{equation}}
\newcommand{\beqn}{\begin{eqnarray}}
\newcommand{\eeqn}{\end{eqnarray}}
\newcommand{\beqns}{\begin{eqnarray*}}
\newcommand{\eeqns}{\end{eqnarray*}}
\newcommand{\bct}{\begin{center}}
\newcommand{\ect}{\end{center}}
\newcommand{\btmz}{\begin{itemize}}
\newcommand{\etmz}{\end{itemize}}
\newcommand{\benum}{\begin{enumerate}}
\newcommand{\eenum}{\end{enumerate}}

\newcommand{\R}{{\mathbb R}}

\newcommand{\lamb}{\bar{\lambda}}

\newcommand{\cL}{\enma{\mathcal L}}

\newcommand{\cN}{\enma{\mathcal N}}
\newcommand{\cR}{\enma{\mathcal R}}

\newcommand{\diag}[1]{ {\sf diag} \! \left( #1 \right) \rule{0em}{1em}}

\newcommand{\lambdab}{{\bar{\lambda}}}

\newcommand{\obtd}[2]{\matbegin \begin{array}{c:c}
        	#1 & #2    \end{array} \matend }	

\newcommand{\tbtd}[4]{\matbegin \begin{array}{c:c}
        	#1 & #2 \\ \hdashline  #3 & #4  \end{array} \matend }

\newcommand{\matbegin}{
    \left[
}
\newcommand{\matend}{
    \right]
}

\newcommand{\bbm}{\begin{bmatrix}} 
\newcommand{\ebm}{\end{bmatrix}} 

\newcommand{\bsm}{\left[ \begin{smallmatrix}} 
\newcommand{\esm}{\end{smallmatrix} \right]} 

\newcommand{\bbNm}{\begin{bNiceMatrix}} 				
\newcommand{\ebNm}{\end{bNiceMatrix}} 
\newcommand{\bNA}[1]{ \left[ \begin{NiceArray}{#1} } 		
\newcommand{\eNA}{ \end{NiceArray} \right] }

\newcommand{\matparbegin}{
    \renewcommand{\baselinestretch}{1}
    \renewcommand{\arraystretch}{.5}
    \setlength{\arraycolsep}{.25em}
    \left[
}
\newcommand{\matparend}{
    \right]
}

\newcommand{\mattightbegin}{
    \renewcommand{\baselinestretch}{1}
    \renewcommand{\arraystretch}{.5}
    \setlength{\arraycolsep}{0.15em}
    \left[
}
\newcommand{\mattightend}{
    \right]
}

\newcommand{\obthdtall}[3]{
 	\matparbegin \begin{array}{c:c:c}
        	 \rule{0em}{1em}  &   & \rule{0em}{0em}  \\  
               	#1 & #2 & #3  \\  
        	  \rule{0em}{1em}  & & \rule{0em}{0em}
	\end{array}\matparend}

\newcommand{\thbotall}[1]{
	  \mattightbegin \begin{array}{c}
       	 \rule{0em}{.7em} \\  #1 \\  \rule{0em}{.7em}
       	\end{array} \mattightend }

\newcommand{\obthwide}[1]{
  	\matparbegin \begin{array}{ccc}
       	 & #1 & 
       	\end{array} \matparend }

\newcommand{\thbodwide}[3]{
  	\matparbegin \begin{array}{ccc}
	 & #1 &     {\rule{0em}{.1em}}_{\rule{0em}{.7em}} 	\\ 	\hdashline
       	\rule{1em}{0em}     & #2 & 		\rule{1em}{0em}	\\ 	\hdashline
	 & #3 &   \rule{0em}{0.9em}
       	\end{array} \matparend }

\newcommand{\thbthd}[9]{
 	\matbegin \begin{array}{c:c:c}
        	#1 & #2 & #3     {\rule{0em}{.1em}}_{\rule{0em}{.7em}}   \\  \hdashline
                #4 & #5 & #6 \\  \hdashline
        	#7 & #8 & #9 \rule{0em}{1em}
	\end{array}\matend}

\newcommand{\lb}{\left(}
\newcommand{\rb}{\right)}

\newcommand{\be}{\begin{equation}}
\newcommand{\ee}{\end{equation}}

\newcommand{\cplxs}{ C\kern -.35em \rule{0.03 em}{.7 ex}~   }

\def\complex{\hbox{C\kern -.45em \rule{0.03 em}{1.5 ex}}~}

\newcommand{\bi}{\begin{itemize}}
\newcommand{\ei}{\end{itemize}}
\newcommand{\ben}{\begin{enumerate}}
\newcommand{\een}{\end{enumerate}}

\newcommand{\tr}{\mbox{tr}}

\newcommand{\inprod}[2]{\left< #1 , #2 \right>}

\newcommand{\bseq}{\begin{subequations}}
\newcommand{\eseq}{\end{subequations}}

\newcommand{\ba}{\begin{array}}
\newcommand{\ea}{\end{array}}

\newcommand{\tcr}[1]{\textcolor{red}{#1}}

\def\clap#1{\hbox to 0pt{\hss#1\hss}}

\newcommand{\btc}{\begin{tabular}{c}}
\newcommand{\btbl}{\begin{tabular}{l}}
\newcommand{\et}{\end{tabular}}

\newcommand{\rem}{\rule{0em}{1em}}
\newcommand{\Diag}[1]{ {\sf dg} \! \left( #1 \right) \rule{0em}{1em}} 
\renewcommand{\diag}[1]{ {\sf diag} \! \left( #1 \right) \rule{0em}{1em}} 
\renewcommand{\tr}[1]{{\sf tr}  \!\left( #1 \rule{0em}{.88em}  \right)} 
\newcommand{\BPi}{{\mathbf \Pi}}
\newcommand{\Bpp}{{\mathbf \Pi}^{\dagger\circ}}

\newcommand{\Lamt}{\tilde{\Lambda}}

\title{A Tutorial on Matrix Perturbation Theory \\ (using compact matrix notation)} 

\author{Bassam Bamieh\thanks{Department of Mechanical Engineering, University of California at Santa Barbara, {\em bamieh@ucsb.edu.} This work is partially supported by NSF Awards CMMI-1763064 and ECCS-1932777.}}

\date{}

\begin{document}

\maketitle

\begin{abstract}
Analytic perturbation theory for matrices and operators is an immensely useful mathematical technique. Most elementary introductions
to this method have their background in the physics literature, and quantum mechanics in particular. In this note, we give an 
introduction  to this method that is independent of any physics notions, and relies purely on concepts from linear algebra. An additional 
feature of this presentation is that matrix notation and methods are used throughout. In particular, we formulate the equations for each term in the analytic 
expansions of eigenvalues and eigenvectors as {\em matrix equations}, namely Sylvester equations in particular. Solvability conditions and 
explicit expressions for solutions of such matrix equations are given, and expressions for each term in the analytic expansions are given 
in terms of those solutions. This unified treatment simplifies somewhat the complex notation that is commonly seen in the literature, and in 
particular, provides relatively compact expressions for the non-Hermitian and  degenerate cases, as well as for higher order terms. 

\end{abstract}

\section{Introduction} 

We want to study the behavior of eigenvalues and eigenvectors of matrices that are a function of a ``small'' parameter $\epsilon$ of the form 
\be
	A_\epsilon ~=~ A_{0} ~+~ \epsilon A_{1}  , 
   \label{Ae.eq}
\ee
where $A_0$ and $A_1$ are given matrices. If the eigenvectors and eigenvalues are analytic functions of $\epsilon$ in a neighborhood of zero, then we can write the eigenvalue/eigenvector relations as power series in $\epsilon$, equate terms of same powers in $\epsilon$ and derive expressions for each set of terms. These  expressions and the related algebra can get messy. It is shown in this document that by adopting matrix notation, expressions are simplified and compactified, and additional insight is obtained. It is rather ironic that {\em matrix notation} is not fully utilized in standard treatments of {\em matrix} perturbation theory. It can be argued that a better way to think about finding expansion terms in the eigenvectors is to treat them all together as a matrix, 
rather than as individual vectors. Better insight is achieved in this manner, especially for the cases of degenerate eigenvalues and higher order terms. 

\subsection*{Notation and Preliminaries} 

Before we begin, we set some useful matrix notation for  eigenvector/eigenvalue relations, as well as manipulations with diagonal 
	matrices and the Hadamard product. 

\begin{itemize}

	\item 
		A vector $v_i$ ($w_i^*$) is a right (left) eigenvector of a matrix $A$ if 
		\[
			A v_i ~=~ \lambda_i v_i, ~~~~~~~\lb w_i^* A ~=~ \lambda w_i^* \rb . 
		\]
		Throughout this note, we will assume the semi-simple case, i.e. that $A_\epsilon$ has a full set of eigenvectors (i.e. 
		diagonalizable) for each $\epsilon$ in some neighborhood of zero. In this case, for an $n\times n$ matrix $A$, 
		 there are $n$ eigenvalue/vector relations
                    \[
                    	Av_i ~=~ \lambda_i v_i 
			\hspace{5em} \lb w_i^* A ~=~ \lambda_i w_i^* \rb , 
			\hspace{5em} i=1,\ldots, n, 
                    \]
                    It is very useful to note that these $n$ relations can be  compactly rewritten as a matrix equation
            \begin{align}
            		 \obthdtall{Av_1}{\cdots}{Av_n}
            		 &=  
            		\obthdtall{\lambda_1 v_1}{\cdots}{\lambda_n v_n} 					\nonumber		\\
            		\Leftrightarrow~~
            		\bbm  & & \\ & A & \\ & & \ebm \obthdtall{v_1}{\cdots}{v_n} 
            		& =  
            		 \obthdtall{v_1}{\cdots}{v_n}
            		\thbthd{\lambda_1}{}{}{}{\ddots}{}{}{}{\lambda_n}					   	
			~~~\Leftrightarrow~~~		
			AV  =  V \Lambda, 												\label{AV.eq}	
            \end{align}
            	where $V$ is a matrix whose columns are the eigenvectors of $A$, and $\Lambda$ 
            	is the diagonal matrix made up of the eigenvalues of $A$. 
		Similarly, we also have 
		\[
			W^* A ~=~ \Lambda W^* , 
		\]
		where the rows of $W^*$ are the left eigenvectors. Note that in the special case of Hermitian (or more generally 
		normal) matrices, the left and right eigenvectors coincide, i.e. $W=V$. 
            		
		\item 	

                    For any square matrix $M$, let $\Diag{M}$ be itself a square diagonal matrix made up of  only the diagonal elements of $M$. Clearly $\Diag{M+N}=\Diag{M}+\Diag{N}$ for any two matrices $M$ and $N$.  If $D$ is a {\em diagonal} matrix with compatible dimensions, then it immediately follows that $\Diag{MD}=\Diag{DM}=\Diag{M}\Diag{D}$. Thus for any two matrices $M$ and $N$ with compatible dimensions 
                    \[
                    	\Diag{M ~\Diag{N}} ~=~ \Diag{M} \Diag{N}.
                    \]

		\item 
		
                        The Hadamard product $M\circ N$ of two matrices is the element-by-element product. It is distributive over matrix additions, but not over matrix products generally except with diagonal matrices
                        \[
                        	M\circ\big(N_1+N_2) ~=~ M\circ N_1+ M\circ N_2 , 
                        	~~~~~~~
                        	M\circ\big(N\Lambda) ~=~ \big( M\circ N\big) \Lambda, 
                        \]
                        when $\Lambda$ is diagonal. If $M=uw^*$ is rank 1, then it follows that 
                        \[
                        	M\circ N ~=~ (uw^*)\circ N ~=~ \diag{u} N~ \diag{w}, 
                        \]
                        where $\diag{u}$ is the diagonal matrix formed from the entries of the vector $u$. Thus if $M=\sum_i \lambda_i u_i w_i^*$ is a diagonalizable 
                        matrix, then we can obtain an expression for the Hadamard product in terms of a sum of standard matrix products
                        \[
                        	M\circ N ~=~ \left( \sum_i \lambda_i u_iw_i^* \right)\circ N  ~=~ \sum_i \lambda_i ~\diag{u_i} N ~\diag{w_i} 
                        \]
                        
                        We will use the notation $M^{k\circ}$ for the Hadamard $k$ power of $M$, i.e. the matrix whose entries are $m_{ij}^k$, with $m_{ij}$
                        being the entries of $M$. Similarly, $M^{-\circ}$ is the matrix with entries $1/m_{ij}$ if all entries are $m_{ij}\neq 0$. 
                        The Hadamard pseudo-inverse 
                        $M^{\dagger\circ}$ is the matrix with entries $1/m_{ij}$ if $m_{ij}\neq 0$, and $0$ if $m_{ij}=0$. 
\end{itemize}                        

\section{Perturbation Expansion} 

Consider~\req{Ae} and   assume that the eigenvectors and eigenvalues are analytic functions of $\epsilon$ in a neighborhood of $0$. The eigenvector/eigenvalue relationship then reads as 
\begin{align}
	A_\epsilon ~V_\epsilon ~=~ V_\epsilon ~ \Lambda_\epsilon, ~~& ~~
	W^*_\epsilon~A_\epsilon  ~=~  \Lambda_\epsilon    W^*_\epsilon  				\nonumber 		\\ 
											& \Updownarrow 											\nonumber		\\ 
	\left( A_{0} ~+~ \epsilon A_{1}   \right)
	\left( V_{0} ~+~ \epsilon V_{1}  ~+~ \cdots  \right)
											& =~ 
	\left( V_{0} ~+~ \epsilon V_{1}  ~+~ \cdots  \right)
	\left( \Lambda_{0} ~+~ \epsilon \Lambda_{1}  ~+~ \cdots  \right)			\label{pertexp.eq}		\\
	\left( W^*_{0} ~+~ \epsilon W^*_{1}  ~+~ \cdots  \right)
	\left( A_{0} ~+~ \epsilon A_{1}   \right)
											& =~ 
	\left( \Lambda_{0} ~+~ \epsilon \Lambda_{1}  ~+~ \cdots  \right)		
	\left( W^*_{0} ~+~ \epsilon W^*_{1}  ~+~ \cdots  \right)									\label{pertexpW.eq}	
\end{align}
Note that $\Lambda$ and $\Lambda_i$'s are {\em diagonal matrices}. 
These equations describe eigenvectors that each belong to a one dimensional subspace, and thus $V_\epsilon$ and $W_\epsilon$ are not unique unless we impose some normalization constraint. There are many possible such constraints. The first we will use is the reciprocal basis constraint $W^*(\epsilon)V(\epsilon)=I$ which gives 
\be
			W^*(\epsilon)V(\epsilon)=I
			~~~~~\Leftrightarrow~~~~~
			\left( W^*_{0} ~+~ \epsilon W^*_{1}  ~+~ \cdots  \right)  \left( V_{0} ~+~ \epsilon V_{1}  ~+~ \cdots  \right)   ~=~I. 
	\label{normexp.eq}		
\ee

Equating equal powers of $\epsilon$ in~\req{pertexp} gives a sequence of matrix equations
\begin{align} 
	A_0 V_0 &=~ V_0\Lambda_0			&							
					W^*_0A_0  &=~ \Lambda_0  W^*																										\label{Lam0.eq}		\\ 
	A_0 V_1 +A_1V_0 &=~ V_0\Lambda_1 + V_1\Lambda_0 	 & 						
					W_1^*A_0  + W_0^* A_1 &=~ \Lambda_1 W_0^* + \Lambda_0 W_1^* 												\label{Lam1.eq}		\\ 
	A_0 V_2 + A_1V_1  &=~  V_0\Lambda_2 + V_1\Lambda_1+ V_2\Lambda_0	 &
					W_2^* A_0  + W_1^* A_1  &=~  \Lambda_2  W_0^* + \Lambda_1 W_1^*+ \Lambda_0 W_2^*			\label{Lam2.eq}		\\ 
								&\vdots																							\nonumber				\\
	A_0 V_{k} +A_1 V_{k-1} &= ~  \sum_{i=0}^k V_i \Lambda_{k-i}	 &
					W_k^* A_0  + W_{k-1}^* A_1  &= ~  \sum_{i=0}^k  \Lambda_{k-i}  W^*_i										\label{LamN.eq}	
\end{align} 
A similar exercise for~\req{normexp} gives 
\be
		W_0^*V_0=I, ~~~~~ W_0^*V_1+W_1^*V_0 = 0, ~~~\cdots~~~  \sum_{i=0}^k W^*_iV_{k-i}=0.
   \label{normcond.eq}
\ee
In the special case of Hermitian matrices, we have $W_\epsilon=V_\epsilon$, and the eigenvectors are thus normalized to be 
orthonormal. We also then have the condition $V_0^*V_1=-V_1^*V_0$, i.e. $V_0^*V_1$ is skew-Hermitian, which implies that its 
diagonal entries must be imaginary, or zero if the vectors are real. In the real case, the diagonal entries are $v_{0i}^*v_{1i}=0$, i.e. for 
each $i$, $v_{0i}$ and $v_{1i}$ are orthogonal.  Geometrically, this means that the curve $v_i(\epsilon)$ has a tangent at $\epsilon=0$ that
is orthogonal to the direction of the unperturbed eigenvector $v_{0i}$. More generally, the normalization condition insures that the curves 
$v_i(\epsilon)$ lie on the sphere in $\R^n$ for any $\epsilon$ (thus the orthogonality of their initial tangents to each $v_{0i}$). 

\section{Calculating First order  Terms} 

	We first begin by calculating the first order behavior of the eigenvalues. This is the term $\Lambda_1$, for which we can use 
	equation~\req{Lam1}. Left multiplication by $W^*_0$ (to get rid of the $V_0$ factor
    multiplying $\Lambda_1$) gives  
    \begin{align*} 
    	W^*_0 A_0 V_1 + 	W^*_0 A_1 V_0 &=~  W_0^* V_0 \Lambda_1  		+ W_0^* V_1 \Lambda_0					\\ 
	\Leftrightarrow~~~~~
    		\Lambda_0  W_0^* V_1 		- W_0^* V_1 \Lambda_0		+	W^*_0 A_1 V_0	&=~   	\Lambda_1    .
    \end{align*} 	
    Now make the following observations 
    \begin{itemize} 
    	\item Since $\Lambda_0$ is diagonal, then $(W_0^*V_1)~\Lambda_0$ and $\Lambda_0~(W_0^*V_1)$ have equal diagonals, i.e. $W_0^*V_1\Lambda_0  - \Lambda_0 W_0^*V_1$ has zeros on the diagonal. 
    	\item Since $\Lambda_1$ must be diagonal, then 
    			\begin{align}
    				\Lambda_1  &=~  \Diag{\rule{0em}{1em} W_0^*A_1V_0  ~-~ W_0^*V_1\Lambda_0  + \Lambda_0 W_0^*V_1} 
    									~=~   \Diag{\rule{0em}{1em} W_0^*A_1V_0 } + 0 					\nonumber		\\
    				\Rightarrow ~~~ \Aboxed{	\Lambda_1	&=~  \Diag{\rule{0em}{1em} W_0^*A_1V_0 } }
				\label{Lam_o.eq}
    			\end{align}
    \end{itemize} 
    Thus the first order correction to the $i$'th eigenvalue is $\lambda_{1i} = w_{0i}^* A_1 v_{0i}$, which is 
    the well known expression\cite{baumgartel1985analytic}. To fully appreciate this, expand~\req{Lam_o} using partitioned matrix notation 
    \[
    	\Lambda_1  ~=~ \Diag{  \thbodwide{w^*_{01}}{\vdots}{w^*_{0n}} 				
				\bbm  & & \\ & A & \\ & & \ebm	\obthdtall{v_{01}}{\cdots}{v_{0n}} }
			~=~ \bbm \lb w_{01}^* A_1 v_{01}\rb  & & \\ & \ddots &  \\  & & \lb w_{0n}^* A_1 v_{0n}\rb \ebm   .
    \]

\subsection{Calculating First Order Eigenvector Terms:  Distinct Eigenvalues} 

Equations~\req{Lam1} can be rewritten as  matrix equations with $V_1$ and $W_1$  as the unknowns respectively  as follows 
\begin{align}
        	A_0 V_1 +A_1V_0 ~=~ V_0\Lambda_1 + V_1\Lambda_0 
	~~~~~~\Leftrightarrow~~~~~~ 
        	A_0 \tcr{V_1}  - \tcr{V_1}\Lambda_0  &=~  \left( V_0\Lambda_1  -A_1V_0	 \right)		,	\label{Sylvester.eq}	\\
		\tcr{W_1^*} A_0   -     \Lambda_0 \tcr{W^*_1} &=~  \left( \Lambda_1 W_0^*  - W_0^* A_1	 \right)			\label{Sylvester2.eq}
\end{align}
The last two equations are Sylvester equations for the matrices $V_1$ and $W_1$ respectively.  We first discuss~\req{Sylvester}. Define the 
matrix-valued operator $\cL(X) ~:=~ A_0X-X\Lambda_0$, and note that~\req{Sylvester} can be rewritten as
                      $
                      		\cL(V_1)  ~=~ \left( V_0\Lambda_1  -A_1V_0	 \right)   .
                      $
The properties of this ``Sylvester operator''   are discussed in Appendix~\ref{Sylvester.app}, and the 
solvability of  equation~\req{Sylvester} is determined by those properties. The following matrices appear in the 
solution: $V_0$, $W_0$, the eigenvectors of $A_0$, $A_0^*$ respectively, and the  matrices 
whose $ij$'th entries are given by 
\[
	\big( \BPi\big)_{ij} ~:=  {\lambda_{0i} - \lambda_{0j} } , 
	~~~~~~~~~~~~
	\left( \BPi^{\dagger\circ} \right)_{ij} ~:= 
		\left\{ 	\begin{array}{ll}		\displaystyle	\frac{1}{\lambda_{0i} - \lambda_{0j} } , ~~~~ & i\neq j, 			\\ 
														0															, 			& i=j. 
					\end{array}			\right.
\]
Entries of $\BPi$ are made up of all possible sums of the eigenvalues of $A_0$ and $-\Lambda_0$ (i.e. differences of the eigenvalues of $A_0$), 
and $\BPi^{\dagger\circ}$ is the Hadamard (element-by-element) pseudo-inverse of $\BPi$.   We now apply the machinery in the Appendices 
to the particular operator $\cL$ in equation~\req{Sylvester}.        
 \begin{itemize} 
 	\item The operator $\cL$ has the following spectral decomposition (applying formula~\req{specform} with $U=Z=I$) 
			\be
				\cL(X) ~=~ 
				V_0  \left(  \BPi \circ \big( W_0^* X \big)   \rule{0em}{1em} \right) 
				~=~  \sum_{ij} (\lambda_{0i} - \lambda_{0j})~  {\big( w_{0i}^*Xe_{0j} \big)} ~v_{0i} e^*_j .
				\label{cLspecdecom.eq}
			\ee
			Since for $i=j$, $\lambda_{0i} - \lambda_{0j}=0$, $\cL$ is rank deficient by at least $n$. 
 	\item Even though $\cL$ is not of full rank,  equation~\req{Sylvester} is solvable since its right hand side is in the 
			range of $\cL$. This argument is detailed in Appendix~\ref{SylEq.app}.
	\item The minimum norm solution of~\req{Sylvester} is obtained from the formula (eq.~\req{cLpinv}, where we again use $U=Z=I$)  
				for the pseudo-inverse $\cL^{\dagger}$ 
				\begin{align}
					V_1 &=~ \cL^{\dagger}  \left( V_0\Lambda_1  -A_1V_0	 \right)	
					~=~  V_0  \left(  \BPi^{\dagger\circ} \circ \big( W_0^*   \left( V_0\Lambda_1  -A_1V_0	 \right)	     \big)   \rule{0em}{1em} \right) 
					~=~  V_0  \left(  \BPi^{\dagger\circ} \circ \big(     \Lambda_1  -W_0^*A_1V_0		     \big)   \rule{0em}{1em} \right) 			\nonumber	\\
					\hspace*{-2em}
						\Rightarrow ~ \Aboxed{ V_1 &=~ -V_0  \left(  \BPi^{\dagger\circ} \circ \big(   W_0^*A_1V_0		 \big)   \rule{0em}{1em} \right) }~ ,	
				\label{SylSol.eq}
				\end{align}
			where we used $W_0^*V_0=I$, the distributive property of the Hadamard product $\circ$, 
			and the fact that $\BPi^{\dagger\circ}\Lambda_1=0$
			since $\Lambda_1$ is diagonal and $\BPi^{\dagger\circ}$ has zeros on the diagonal. 
		\item 
			The Sylvester operator in~\req{Sylvester2} is just $\cL^*$, and we can similarly derive (see Appendix~\ref{SylEq.app})
			 the minimum-norm 
			solution
				\[
					\boxed{ W_1^* =~   \left(  \BPi^{\dagger\circ} \circ \big(   W_0^*A_1V_0		 \big)   \rule{0em}{1em} \right) W^*_0}~ .
				\]	

\end{itemize}

To compare the  solution formula~\req{SylSol} with standard expressions in the literature, we expand it as 
\[
	V_1 ~=~- V_0 \left(  \BPi^{\dagger\circ} \circ \big(   W_0^*A_1V_0		     \big)   \rule{0em}{1em} \right) 
			~=~ -\sum_{i\neq j} \frac{1}{\lambda_{0i} - \lambda_{0j}}  {\big( w_{0i}^*A_1 v_{0j} \big)} 
					~v_{0i} e^*_j 
\]
Note that $V_1$ is a matrix whose $k$'th column is $v_{1k}$. The $k$'th column 
is obtained from $V_1e_k$ 
\be
		v_{1k} ~=~- V_1e_k ~=~ -\sum_{i\neq j}   \frac{w_{0i}^*A_1v_{0j}}{\lambda_{0i}-\lambda_{0j}}  ~v_{0i} e_j^* ~e_k
					~=~  -\sum_{i\neq k}   \frac{w_{0i}^*A_1v_{0k}}{\lambda_{0i}-\lambda_{0k}}  ~v_{0i}, 
					~=~  \sum_{i\neq k}   \frac{w_{0i}^*A_1v_{0k}}{\lambda_{0k}-\lambda_{0i}}  ~v_{0i}, 
\ee
which is the standard expression in the literature~\cite{baumgartel1985analytic}. 


\subsubsection*{Uniqueness of $V_1$ and $W_1$}

Note that~\req{SylSol} is just one solution to the Sylvester equation. All other solutions can be obtained by adding 
arbitrary elements of the null space of $\cL$. In  Appendix~\ref{SylEq.app},  it is  shown how the null spaces of 
			$\cL$ and $\cL^*$ are characterized, and from that  we can conclude that 
all solutions of the Sylvester equation~\req{Sylvester} and~\req{Sylvester2} can be written as
\be
	 V_1 ~=~  -V_0 \left( \BPi^{\dagger\circ} \circ \big( W_0^* A_1 V_0 \big) \rule{0em}{1em} \right)   + V_0D_1, 
	 ~~~~~
	 W_1^* ~=~   \left(  \BPi^{\dagger\circ} \circ \big(   W_0^*A_1V_0		 \big)   \rule{0em}{1em} \right) W^*_0 -D_2W_0^*,
  \label{allsol.eq}	 
\ee
where $D_1$ and $D_2$ are arbitrary diagonal matrices. 
The normalization conditions~\req{normcond} now give  
\begin{align}
	 0 ~=~ W_0^*V_1+W_1^*V_0 
	 	& =~ - W_0^*V_0 \left( \BPi^{\dagger\circ} \circ \big( W_0^* A_1 V_0 \big) \rule{0em}{1em} \right)   +
	 				 \left( \BPi^{\dagger\circ } \circ \big( W_0^* A_1 V_0 \big) \rule{0em}{1em} \right) W_0^*V_0  + D_1-D_2 		\nonumber	\\
	 	& =~ 			D_1-D_2 .																																						\label{d1d2.eq}
\end{align}
We now make several observations
\begin{itemize} 
	\item 
			In the self-adjoint case, we have $W_0=V_0$ and $W_1=V_1$, and the normalization conditions~\req{normcond} give
			\[
				0 ~=~ V_0^*V_1 
				~=~ V_0^* 	\left( -V_0 \left( \BPi^{\dagger\circ} \circ \big( W_0^* A_1 V_0 \big) \rule{0em}{1em} \right)   + V_0D_1 \right)
				~=~ -\BPi^{\dagger\circ} \circ \big( W_0^* A_1 V_0 \big) \rule{0em}{1em}   + D_1  .
			\]
			Since the first term has zeros on the diagonal ($\BPi^{\dagger\circ}$ is zero on the diagonal)  
			and $D_1$ is diagonal, we conclude that $D_1=0$. Thus, the minimum 
			norm solution~\req{SylSol}  to the Sylvester equation does indeed satisfy the normalization conditions~\req{normcond}. 
			
		\item
			In the general (non self-adjoint) case, it seems that the following set of solutions 
              \[
              	 V_1 ~=~  -V_0 \left( \BPi^{\dagger\circ} \circ \big( W_0^* A_1 V_0 \big) \rule{0em}{1em} \right)   + V_0D, 
              	 ~~~~~
              	 W_1^* ~=~   \left(  \BPi^{\dagger\circ} \circ \big(   W_0^*A_1V_0		 \big)   \rule{0em}{1em} \right) W^*_0 -DW_0^*
              \]
              where $D$ is any diagonal matrix will all satisfy the normalization condition~\req{normcond} up to first order. Since 
              $V_1$ and $W_1$ must be 
              unique\footnote{The condition $W_\epsilon^* V_\epsilon=I$ uniquely determine the functions $V_\epsilon$ and $W_\epsilon$, 
              and therefore uniquely determines their Taylor series expansion terms $\{V_k\}$ and $\{W_k\}$.  }, 
              it seems like consideration of higher order terms is necessary to find the 
              correctly normalized solutions. 
	\item 
			In some references~\cite[Eq. 5.1.31]{sakurai1995modern}, a non-standard normalization is used. In our notation, this 
			normalization is stated as follows
			\begin{align}
				w_{0i}^*v_{0j} & =~ \delta_{i-j},    ~~~~~~~~~~ w_{0i}^*v_{\epsilon i} ~=~1, 								\label{nonsnorm.eq}		\\ 
				W_0^*V_0 &=I,    ~~~~~ ~~~~\Diag{W_0^*V_\epsilon} ~=~ I ~~~~~~\Rightarrow~~~ \forall k\geq 1, ~\Diag{W_0^*V_k} ~=~ 0.
				\nonumber
			\end{align}
			This normalization will uniquely determine $v_{\epsilon i}$ unless it becomes orthogonal to $w_{0i}$. For finite matrices, 
			we can guarantee that $v_{\epsilon i}$ will not be orthogonal to $w_{0i}$ for $\epsilon$ in some neighborhood of $0$. 
			
			This non-standard normalization is used because it simplifies the recursive formulae for higher order peturbation 
			terms as we will see later on. For now, observe that if we adopt this normalization and check the general solutions~\req{allsol}
			\begin{align*}
				0 &=~ \Diag{W_0^*V_1} 
				~=~ 
				\Diag{ -W_0^*V_0 \left( \BPi^{\dagger\circ} \circ \big( W_0^* A_1 V_0 \big) \rule{0em}{1em} \right)   + W_0^*V_0D_1  }		\\
				&=~ 
				-\Diag{  \BPi^{\dagger\circ} \circ \big( W_0^* A_1 V_0 \big) }  +  D_1   ~=~ 0~+~D_1 
			\end{align*}
			So we can conclude that the minimum-norm solution~\req{SylSol} is indeed normalized with the 
			non-standard normalization~\req{nonsnorm}. 
              
\end{itemize}

%


\subsection{Calculating First Order Eigenvector Terms $V_1$:  Repeated (Degenerate) Case} 

We consider the case when the eigenvalues are repeated, but $A_0$ has a full set of eigenvectors. This condition is 
automatically satisfied in the Hermitian or normal $A_0$ case, but may not be always true in the non-normal $A_0$ case. We will 
consider the case of normal $A_0,A_1$ here, for which we have $W_0=V_0$, and therefore $A_0=V_0\Lambda_0V^*_0$.

The main difficulty when $A_0$ has repeated
eigenvalues is that there is not a unique choice of the eigenvectors $V_0$ in this case, even after normalization. For any 
repeated eigenvalue, there corresponds an $m$-dimensional invariant subspace, where $m$ is the geometric multiplicity 
of the eigenvalue. Any basis of this invariant subspace is composed of eigenvectors. It turns out that there is a special 
choice of basis which renders the expansion~\req{pertexp} valid. The main task is to find such a basis. 

Now assume $A_0$ has a repeated eigenvalue $\lamb$ (multiplicity $m$) and define a basis so that $A_0$ is block-diagonal 
with the first $m$ basis elements  a basis for the eigensubspace of $\lamb$. With this basis,  $A_0$, 
$V_0$, $W_0$ and $\Lambda_0$ have a block partitioning as 
\[
	A_0 = \bbm  (A_{0})_{11} & 0\\ 0 & (A_{0})_{22} \ebm 
	= \bbm  (V_{0})_{11} & 0\\ 0 & (V_{0})_{22} \ebm 
		\bbm  \lamb I & 0\\ 0 & (\Lambda_{0})_{22} \ebm
		\bbm  (V^*_{0})_{11} & 0\\ 0 & (V^*_{0})_{22} \ebm
	= V_0\Lambda_0 V_0^*, 
\]
where $I$ is the $m\times m$ identity matrix. 
Clearly there is not a unique choice for $(V_0)_{11}$ since replacing it with  $(V_0)_{11}U$, where $U$ is any unitary matrix will keep the above expression valid\footnote{$U$ unitary means $UU^*=I$, and therefore $(V_0)_{11}U~ \lambdab I ~U^* (V_0^*)_{11} = (V_0)_{11} ~\lambdab I ~ (V_0^*)_{11} $.}. 

Now take eq.~\req{Lam1}
\[
	\Lambda_0 V_1 + A_1V_0 ~=~  V_0\Lambda_1 + V_1 \Lambda_0.
\]
If we partition all matrices conformably with the above partitions,  we get 
\begin{multline*}
	\bbm  \lamb I & 0\\ 0 & (\Lambda_{0})_{22} \ebm   
	\bbm (V_1)_{11} & (V_1)_{12} \\ (V_1)_{21}  & (V_1)_{22} \ebm
	+ \bbm (A_1)_{11} & (A_1)_{12} \\ (A_1)_{21}  & (A_1)_{22} \ebm
	\bbm  (V_{0})_{11} & 0\\ 0 & (V_{0})_{22} \ebm																\\
	=~ 
	\bbm  (V_{0})_{11} & 0\\ 0 & (V_{0})_{22} \ebm
	 \bbm (\Lambda_1)_{11} & 0 \\ 0  & (\Lambda_1)_{22} \ebm
	+ 
	\bbm (V_1)_{11} & (V_1)_{12} \\ (V_1)_{21}  & (V_1)_{22} \ebm
	\bbm  \lamb I & 0\\ 0 & (\Lambda_{0})_{22} \ebm .
\end{multline*} 
Taking the $11$ block equation we get 
\[
	\lamb I_n  ~(V_1)_{11}  ~+~ (A_1)_{11} (V_{0})_{11} ~=~(V_{0})_{11} (\Lambda_1)_{11} ~+~ (V_1)_{11} ~\lamb I_n  ,
\]
but $\lamb I_n$ commutes with any matrix, so we simply get 
\[
	(A_1)_{11} (V_{0})_{11} ~=~(V_{0})_{11} (\Lambda_1)_{11},
\]
from which we can obtain $(\Lambda_1)_{11}$ by left multiplication by $(V_{0}^*)_{11}$
\[
	(\Lambda_1)_{11} ~=~ 
	(V_{0}^*)_{11} (A_1)_{11} (V_{0})_{11}  ~=~ (V^*_{0} A_1 V_{0})_{11} , 
\]
where the last equality follows from the block-diagonal structure of $V_0$ and $V_0^*$. 
Now let's examine the meaning of this. The perturbation expansion~\req{pertexp} assumes that all the $\Lambda_i$'s are
diagonal matrices (otherwise~\req{pertexp} is not an eigenvalue/eigenvector relation). This equation states that for $(\Lambda_1)_{11}$ to be diagonal as required, {\em we must choose $(V_0)_{11}$ so that its columns are the eigenvectors 
of $(A_1)_{11}$}. Consequently, $(\Lambda_1)_{11}$ will simply be the diagonal matrix of eigenvalues of $(A_1)_{11}$. 

The meaning of the above is that only this particular choice of $V_0$ will yield the expansion~\req{pertexp}. Unlike the
non-repeated case, we also have to discover ``the special'' $V_0$, not just $V_1$ and $\Lambda_1$.

Now to calculate $V_1$ we can apply the psuedo-inverse formula~\req{cLpinv} (which produced~\req{SylSol} earlier), but now we pay special attention 
to the structure of $\BPi^{\dagger\circ}$. In this case, it is 
\[
	\left( \BPi^{\dagger\circ} \right)_{ij} ~:= 
		\left\{ 	\begin{array}{ll}		\displaystyle	\frac{1}{\lambda_{0i} - \lambda_{0j} } , ~~~~ &\lambda_{0i} \neq \lambda_{0j}, 			\\ 
														0															, 			&\lambda_{0i} = \lambda_{0j}. 
					\end{array}			\right.
\]
It is instructive to look at the structure of $\BPi^{\dagger\circ}$, it is of the following form 
\[
	\BPi^{\dagger\circ} ~=~ \bbm \scalebox{2}{0} &    &    \\ 
																	&  0 &   \\ 
																	&     & \ddots    				\ebm ,
\]
where the large $0$ block is $m\times m$, the diagonal is all zeros, and the remaining terms are $1/(\lambda_{0i} - \lambda_{0j})$. 
The expression for $V_1$ is now the same as~\req{SylSol}, but now we also make sure to use the $V_0$ that made $\Lambda_1$ diagonal 
above (recall that the derivation of~\req{SylSol} assumed $\Lambda_1$ to be diagonal) 
\[
	V_1 ~=~ -V_0  \left(  \BPi^{\dagger\circ} \circ \big(   V_0^*A_1V_0		 \big)   \rule{0em}{1em} \right) . 
\]
Again, to compare with existing expressions, find $v_{1k}$ using $V_1e_k$ 
\be
		v_{1k} ~=~- V_1e_k ~=~ -\sum_{\lambda_{0i}\neq \lambda_{0k}}   \frac{v_{0i}^*A_1v_{0j}}{\lambda_{0i}-\lambda_{0j}}  ~v_{0i} e_j^* ~e_k
					~=~  \sum_{\lambda_{0i}\neq \lambda_{0k}}   \frac{v_{0i}^*A_1v_{0k}}{\lambda_{0k}-\lambda_{0i}}  ~v_{0i}, 
\ee

\begin{example} 
	Consider the simplest case of $2\times 2$ matrices 
	\[
		A_0~+~\epsilon A_1 ~:=~ I ~+~ \epsilon M  , 
	\] 
	where $M$ is a normal matrix. 
	$A_0=I$ has two eigenvalues at $1$ with multiplicity $2$, and a 2-dimensional eigenspace (all of $R^2$). Let $M$ have 
	eigenvalue/vector pairs $(\alpha,x)$ and  $(\beta,y)$. We actually know the eigenvalues of $I+\epsilon M$
	explicitly  as a function 
	of $\epsilon$ as (since $I+\epsilon M$ has the same eigenvectors as $M$) 
	\[
		\lambda_{\epsilon 1}  ~=~ 1~+~ \epsilon~ \alpha , 
		\hspace{4em} 
		\lambda_{\epsilon 2}  ~=~ 1~+~ \epsilon~ \beta .		
	\]
	Let's see if this can be replicated using the procedure described above. We need to choose $V_0$ so that its columns 
	are eigenvectors of $M$, so the choice is 
	\[
		V_0 ~=~ \obtd{x\rule{0em}{1em}}{y} 	, 
	\]
	and therefore 
	\[
		V_0^* A_1 V_0 ~=~ \obtd{x\rule{0em}{1em}}{y}^* M \obtd{x\rule{0em}{1em}}{y}  
			~=~ \tbtd{\alpha}{0}{0}{\beta} .
	\]
	Note that the off-diagonal elements are zero since the vectors $x$ and $y$ are orthogonal ($M$ is normal). 
\end{example}

\section{Higher Order Terms} 

The following normalization will lead to rather compact recursive formulae for terms of all orders 
\begin{align}
	w_{0i}^*v_{0j} & =~ \delta_{i-j},    ~~~~~~~~~~ w_{0i}^*v_{\epsilon i} ~=~1, 							\nonumber		\\ 
	W_0^*V_0 &=I,    ~~~~~ ~~~~\Diag{W_0^*V_\epsilon} ~=~ I ~~~~~~\Rightarrow~~~ \forall k\geq 1, ~\Diag{W_0^*V_k} ~=~ 0.
	\label{nonsnorm2.eq}	
\end{align}
We find the eigenvalues first. Decompose $\Lambda_\epsilon = (\Lambda_0 +\epsilon\Lambda_1 + \cdots)=:(\Lambda_0+\Lamt)$, and 
rearrange to obtain an equation for $\Lamt$ 
\begin{align}				
	\left( A_{0} ~+~ \epsilon A_{1}   \right)
	V_\epsilon
											& =~ 
	V_\epsilon
	\left( \Lambda_{0} ~+~  \Lamt  \right)														&&	\label{VL.eq}		\\
	V_\epsilon  \Lamt
											& =~ 
	-V_\epsilon
	 \Lambda_{0} + A_{0}V_\epsilon + \epsilon A_{1}  V_\epsilon							&& \mbox{rearrange terms} \nonumber	\\
	W_0^* V_\epsilon  \Lamt
											& =~ 
	-W_0^*V_\epsilon
	 \Lambda_{0} + W^*_0A_{0}V_\epsilon + \epsilon W^*_0 A_{1}  V_\epsilon		&& \mbox{left multiply by $W_0^*$}\nonumber \\
	\Diag{ W_0^* V_\epsilon  \Lamt }
											& =~ 
	-\Diag{W_0^*V_\epsilon  \Lambda_{0}}
	  + \Diag{ \Lambda_0~ W^*_0 V_\epsilon} +  \Diag{W^*_0 A_{1}  ~\epsilon V_\epsilon}	&& \mbox{take diagonals}\nonumber 	\\
	\Diag{ W_0^* V_\epsilon}  \Lamt 
											& =~ 
	-\Diag{W_0^*V_\epsilon}
	 \Lambda_{0} +  \Lambda_0~ \Diag{W^*_0 V_\epsilon} +										
	  \Diag{W^*_0 A_{1} ~ \epsilon V_\epsilon}								&& \mbox{$\Diag{M\Lambda}=\Diag{M}\Lambda$ if $\Lambda$ diagonal}\nonumber\\
	  \Lamt 
											& =~ 										
	  \Diag{W^*_0 A_{1} ~ \epsilon V_\epsilon}								&& \mbox{$\Diag{W_0^*V_\epsilon}=I$}\nonumber\\
	\big( \epsilon \Lambda_1 + \epsilon^2 \Lambda_2 + \cdots \big) 
											& =~ 
	  \Diag{W^*_0 A_{1} ~ \epsilon \big( V_0+\epsilon V_1 +\cdots \big) \rule{0em}{1em}   }			&&\mbox{expand \& equate}	\nonumber		\\
	\Aboxed{
	 \Lambda_k 
											& =~ 
	  \Diag{W^*_0 A_{1}  V_{k-1} \rule{0em}{1em}   }} 																			\label{Lamk.eq}
\end{align}

The first of these formulae is the familiar $\Lambda_1= \Diag{W^*_0 A_{1}  V_{0}  }$. Successive $\Lambda_k$'s require finding 
the vectors $V_{k-1}$ that are normalized according to~\req{nonsnorm2} (and not the standard normalization). 

To find $V_k$'s, recall and rearrange equation~\req{LamN} into a Sylvester equation for $V_k$ 
\begin{align*}
	A_0 V_{k} +A_1 V_{k-1} &= ~  V_0\Lambda_k + V_1\Lambda_{k-1} \cdots +V_{k-1} \Lambda_1 + V_k\Lambda_0					\\
	A_0 V_{k} - \Lambda_k V_0  &= ~ 	V_0\Lambda_{k}+V_1\Lambda_{k-1}+ \cdots +V_{k-1} \Lambda_1 ~-~ A_1 V_{k-1}
\end{align*}
All solutions to this Sylvester equation are given by the psuedo-inverse formula 
\[
	V_k ~=~ V_0 \left(  \Bpp \circ  W_0^*  \big(  V_0\Lambda_{k}+ V_1\Lambda_{k-1}+ \cdots +V_{k-1} \Lambda_1 ~-~ A_1 V_{k-1}  \big)\rem     \right)
						+V_0 D_k, 
\]
where $D_k$ is any diagonal matrix. Enforcing the normalization $W_0^*V_0=I$, $\Diag{W_0^*V_k}=0, ~k\geq 1$ gives that
$D_0=I$ and $D_k=0$ for $k\geq 1$, and we conclude 
\[
	\boxed{
	V_k ~=~ V_0 \left(  \Bpp \circ  W_0^*  \big(   V_0\Lambda_{k}+V_1\Lambda_{k-1}+ \cdots +V_{k-1} \Lambda_1 ~-~ A_1 V_{k-1}  \big)\rem     \right)
						}~. 
\]
When this formula is combined with~\req{Lamk} for $\Lambda_k$, we see that we can obtain $V_k$ using previous calculations for the terms $V_0,\ldots,V_{k-1}$ and $\Lambda_0,\ldots,\Lambda_{k-1}$. For computations, the following equivalent form may be easier to implement
\[
	V_k ~=~ V_0 \left( 
			\sum_{i=0}^{k-1} \big(  \Bpp \circ  W_0^*    V_i \big) \Lambda_{i-1}   
				~-~ \Bpp \circ \left(  W_0^* A_1 V_{k-1}  \right)      \right). 
\]

\section{Acknowledgement} 

The author would like to acknowledge helpful input from Maurice Filo, Stacy Patterson and Karthik Chickmaglur. 

\newpage 

\appendix 

\section{Appendix}

\subsection{Spectral Decomposition of the Sylvester Operator} 				\label{Sylvester.app} 

Consider the Sylvester matrix equation (in $X$) of the form 
\[
	AX~+~XB~=~Q, 
\]
where $A$, $B$ and $Q$ are given (square) matrices. Solutions and properties of this equation are determined by the 
Sylvester operator $\cL$ 
\[
	\cL(X) ~:=~ AX~+~XB. 
\]	
We will need to find the ``eigen-matrices'' of $\cL$ and its adjoint $\cL^*$. First we calculate $\cL^*$ using the standard inner
product on matrices $\inprod{M}{N}=\tr{M^*N}$. Starting from the $\equiv$ equality 
{\small
\[
	\tr{  X^*(A^*Y+YB^*)}  =
	\tr{  X^*A^*Y+B^*X^*Y}  =
	\tr{  (AX+XB)^* Y}  =
	\inprod{\cL(X)}{Y} \equiv  \inprod{X}{\cL^*(Y)} 
	=\tr{  X^*\cL^*(Y)}
\]			}
we see that 
\[
	\cL^*(Y) ~=~ A^*Y+YB^*. 
\]
Thus  in particular, $\cL$ is self-adjoint only if $A$ and $B$ are self-adjoint\footnote{A simple calculation will show that $\cL$ is 
normal if $A$ and $B$ are normal.}

Now we assume  that  $A$ and $B$ are  diagonalizable (i.e. have a full set of linearly independent eigenvectors each). 
A spectral decomposition of $\cL$ can be obtained from the spectral decompositions of $A$ and $B$. This will then allow for applying arbitrary 
analytic functions on $\cL$ including inversion. 

Denote the  eigenvalues and right/left eigenvectors of  $A$ and $B$ as follows 
\[
	\begin{array}{rclcrcl}
		Av_i & = & \lambda_i v_i	, & & A^*w_i & = & \lambda_i^*	 w_i ,	\\ 
		Bu_j & = & \gamma_j u_j	, & &   B^* z_j & = & \gamma_j^* z_j,
	\end{array} 
	~~~~~\mbox{or in matrix form}~~~~~
	\begin{array}{rclcrcl}
		AV & = & V\Lambda , & & A^*W & = & W \Lambda^*			\\
		BU & = & U\Gamma , & & B^*Z & = & Z\Gamma^*
	\end{array}	.
\]
Note that $W=V^{-*}$ and $Z=U^{-*}$. 
The $n^2$ eigenvalues and ``eigen-matrices'' of $\cL$ and $\cL^*$ are found as follows 
\begin{alignat*}{3}
	\cL(v_iz_j^*) &=~ A~ v_iz_j^* ~+~ v_iz_j^* ~B &&=~ \lambda_i ~v_iz_j^* ~+~ \gamma_j ~v_iz^*_j &&=~ (\lambda_i + \gamma_j) ~v_iz^*_j,		\\
	\cL^*(w_iu^*_j) &=~ A^*~ w_iu^*_j ~+~ w_iu^*_j ~B^* &&=~ \lambda_i^*~ w_iu^*_j ~+~ \gamma_j^* w_iu^*_j  &&=~ (\lambda_i^* + \gamma_j^*) ~w_iu^*_j. 
\end{alignat*}
In other words, eigenvalues of $\cL$ (resp. $\cL^*$) are all $n^2$ possible combinations $\lambda_i + \gamma_j$ (resp. $(\lambda_i+\gamma_j)^*$ )
of eigenvalues of $A$ and $B$, 
and the eigen-matrices are all the corresponding outer products of right/left eigenvectors of $A$ and $B$. 

Using the above, a spectral decomposition of $\cL$ can now be written as follows 
\[
	\cL(X) ~=~ \sum_{ij} (\lambda_i + \gamma_j)  \inprod{w_iu^*_j}{X} ~v_i z^*_j .
\]
This can be rewritten in compact notation. First define the matrix $\BPi$ whose $ij$'th entry is
\[
	 \big( \BPi  \big)_{ij} ~:=~ \lambda_i +\gamma_j. 
\] 
Noting that  $v_i=Ve_i$ (where $e_i$ is the vector 
of all zeros except $1$ in the $i$'th row), and similarly for the other eigenvectors, calculate 
\begin{align}
	\cL(X) &=~ \sum_{ij} (\lambda_i + \gamma_j)  \inprod{w_iu^*_j}{X} ~v_i z^*_j 
				~=~ \sum_{ij} (\lambda_i + \gamma_j)~  {\big( w_i^*Xu_j \big)} ~Ve_i e^*_j Z^*			\nonumber			\\
				&=~ V \left(  \sum_{ij} (\lambda_i + \gamma_j)~  \big( w_i^*Xu_j \big) ~e_i e^*_j   \right)  Z^*
				~=~ V  \left(  \BPi \circ \big( W^* X U \big)   \rule{0em}{1em} \right)  Z^*, 						\label{specform.eq}
\end{align}
where $\circ$ is the Hadamard (element-by-element) product of matrices. The last equation follows from observing that 
$W^*XU$ is the matrix whose $ij$'th entry is $ \big( w_i^*Xu_j \big)$. 
For later reference, we also calculate $\cL^*$ 
\begin{align}
	\cL^*(X) &=~ \sum_{ij} (\lambda^*_i + \gamma^*_j)  \inprod{v_i z^*_j}{X} ~ w_iu^*_j
				~=~ \sum_{ij} (\lambda^*_i + \gamma^*_j)~  {\big( v_i^*Xz_j \big)} ~We_i e^*_j U^*			\nonumber			\\
				&=~ W \left(  \sum_{ij} (\lambda_i + \gamma_j)~  \big( v_i^*Xz_j \big) ~e_i e^*_j   \right)  U^*
				~=~ W  \left( \bar{ \BPi} \circ \big( V^* X Z \big)   \rule{0em}{1em} \right)  U^*, 						\label{specformLs.eq}
\end{align}
where $\bar{\BPi}$ is the complex conjugate (without transposing) of $\BPi$. 


The inverse $\cL^{-1}$ (if it exists) and the pseudo-inverse $\cL^{\dagger}$ can be calculated from the spectral decomposition 
\[
	\cL^{-1}(X) ~=~ \sum_{ij} \frac{1}{\lambda_i + \gamma_j}  \inprod{w_iu^*_j}{X} ~v_i z^*_j , 
	~~~~~~~
	\cL^{\dagger}(X) ~=~ \sum_{\stackrel{ij}{ \lambda_i + \gamma_j \neq 0} } \frac{1}{\lambda_i + \gamma_j}  \inprod{w_iu^*_j}{X} ~v_i z^*_j .	
\]
We can also give compact formulae if we define $\BPi^{-\circ}$ and $\BPi^{\dagger\circ}$ using element-by-element operations 
\[
	\left( \BPi^{-\circ} \right)_{ij} ~:=  \frac{1}{\lambda_{i} + \gamma_{j} } , 
	~~~~~~~~~~~~
	\left( \BPi^{\dagger\circ} \right)_{ij} ~:= 
		\left\{ 	\begin{array}{ll}		\displaystyle	\frac{1}{\lambda_{i} + \gamma_{j} } , ~~~~ & \lambda_{i} + \gamma_{j} \neq 0, 			\\ 
														0															, 			&  \lambda_{i} + \gamma_{j} =0. 
					\end{array}			\right.
\]
$\cL^{-1}$ and $\cL^\dagger$ can now be rewritten as 
\begin{align}
	\cL^{-1}(X) 	&=~ V  \left(  \BPi^{-\circ} \circ \big( W^* X U \big)   \rule{0em}{1em} \right)  Z^*, 			\nonumber		\\ 
	\cL^{\dagger}(X) 	&=~ V  \left(  \BPi^{\dagger\circ} \circ \big( W^* X U \big)   \rule{0em}{1em} \right)  Z^*. 
  \label{cLpinv.eq}	
\end{align}

The above can be generalized using the spectral decomposition to any function $f$ analytic in a neighborhood of the 
spectrum of $\cL$ by 
\[
	\big(f(\cL)\big)(X) 	~=~ V  \left(  f^\circ(\BPi) \circ \big( W^* X U \big)   \rule{0em}{1em} \right)  Z^*, 
\]
where $f^\circ(\BPi)$ is the element-by-element application of the function $f$ on each entry of the matrix $\BPi$. 
For example, given a matrix differential equation 
\[
	\dot{X} ~=~ AX+XB ~=~ \cL(X) , 	~~~~~~~~X(0)=\bar{X}
\]	
we can write the solution formally as $\left(e^{t\cL} \right)(\bar{X})$. The formula above gives 
\[
	X(t) ~=~ \left(e^{t\cL} \right)(\bar{X}) ~=~ V  \left(  e^{\circ(t\BPi)} \circ \big( W^* \bar{X} U \big)   \rule{0em}{1em} \right)  Z^*, 
\]
where  $e^{\circ(t\BPi)}$ is the matrix whose $ij$'th entry is $e^{t(\lambda_i+\gamma_j)}$.

\subsection{Solvability of the Sylvester Equations~\req{Sylvester} and~\req{Sylvester2}}				\label{SylEq.app}

We begin with the equation~\req{Sylvester} for $V_1$. 

		\bi
          	\item We need to characterize $\cN(\cL)$,  the null space of $\cL$. Using~\req{cLspecdecom} we note that $\cL(X)=0$ means 
          			\[
          					 \forall i,j, ~  (\lambda_{0i}-\lambda_{0j}) ~( w_{0i}^* X e_j  ) = 0
					 		~~~~\Longleftrightarrow~~~~ 
							\forall i\neq j, ~( w_{0i}^* X e_j  ) = 0
 					 		~~~~\Longleftrightarrow~~~~ 
							W_0^*X ~=~ D, 
         			\]
					where $D$ is some diagonal matrix (the diagonal entries of $D$ are the numbers $(w_{0i}^* X e_i  )$ which 
					cannot be determined from the above condition). Finally we note that $W_0$ and $V_0$ are inverses of each other, 
					so  $W_0^*X=D$ is equivalent to $X=V_0D$ and we conclude 
%
          			\[
          				X\in\cN(\cL) ~~~~\Leftrightarrow~~~~  X ~=~ V_0 D, ~~\mbox{$D=$ any diagonal matrix}
          			\]
			
			\item Similarly for $\cL^*(X)~=~ A_0^*X-X\Lambda_0^*$. It's spectral decomposition is 
					\[
							 \sum_{ij} (\lambda_{0i}^* - \lambda_{0j}^*)~  {\big( v_{0i}^*Xe_{j} \big)} ~w_{0i} e^*_j, 
					\]	
					and again $\cL^*(X)=0$ iff $\forall i\neq j, ~~{\big( v_i^*Xe_j \big)} =0$.  Repeating the above argument we conclude 
          			\[
          				X\in\cN(\cL^*) ~~~~\Leftrightarrow~~~~  X ~=~ W_0 D, ~~\mbox{$D=$ any diagonal matrix}
          			\]
			
          	
          	\item 
                      The Sylvester equation~\req{Sylvester}  is solvable iff 
                      $
                      	\left( V_0\Lambda_1  -A_1V_0	 \right)~~ \in~~ \cR(\cL)  =  \cN(\cL^*)^{\perp} .
                      $
                      Thus to verify solvability, we need to show that $V_0\Lambda_1  -A_1V_0	$ is orthogonal to $\cN(\cL^*)$. Indeed 
                      \begin{align*}
                      	X\in\cN(\cL^*) &\Leftrightarrow ~ X =W_0D,  \\ 
                      	 \inprod{X}{V_0\Lambda_1  -A_1V_0	} &= \inprod{W_0D}{V_0\Lambda_1  -A_1V_0	}
                      	 		~=~ \tr{D^*W_0^*V_0\Lambda_1}- \tr{D^*W_0^*A_1V_0}			\\
                      			&=~ \tr{D^*\Lambda_1}- \tr{D^*~W_0^*A_1V_0}~=~ \tr{D^*~ \Diag{W_0^*A_1V_0} }- \tr{D^*~W_0^*A_1V_0} \\
								&=~ \tr{D^*~ {W_0^*A_1V_0} }- \tr{D^*~W_0^*A_1V_0} ~=~0,
                      \end{align*}
                      where we used $\Lambda_1=\Diag{W_0^*A_1V_0}$, and that $D$ is a diagonal matrix. 
          		\ei

For equation~\req{Sylvester2}, note that by transposing it we get a Sylvester equation of a similar form to~\req{Sylvester} 
\[
	A_0^*W_1 -W_1\Lambda_0^* ~=~ W_0\Lambda_1^* - A_1^*W_0. 
\]
This equation is $\cL^*(W_1) = W_0\Lambda_1^* - A_1^*W_0$ where $\cL^*$ is the adjoint of the Sylvester operator of equation~\req{Sylvester}, which we have already analyzed. In particular, the spectral decomposition of $\cL^*$ can be calculated from~\req{specformLs}
\[
	\cL^*(Y) ~=~ W_0  \left(  \bar{\BPi} \circ \big( V^*_0 Y  \big)   \rule{0em}{1em} \right) . 
\] 
Applying the pseudo-inverse to the right hand side 
\begin{align*}
		W_1 &=~
		\cL^{*\dagger}(W_0\Lambda_1^* - A_1^*W_0)  
		~=~ W_0  \left(  {\bar{\BPi}}^{\dagger\circ} \circ \big( V_0^* \left(  W_0\Lambda_1^* - A_1^*W_0 \right)    \big)   \rule{0em}{1em} \right) 
		~=~ W_0  \left(  {\bar{\BPi}}^{\dagger\circ} \circ \big(    \Lambda_1^* -V_0^* A_1^*W_0     \big)   \rule{0em}{1em} \right) \\ 
		&=~ -W_0  \left(  {\bar{\BPi}}^{\dagger\circ} \circ \big(    V_0^* A_1^*W_0     \big)   \rule{0em}{1em} \right) ,
\end{align*}
where the last equality follows from $\bar{\BPi}^{\dagger\circ}$ having all zeros on the diagonal. 
We can rewrite the solution for $W_1^*$ by noting that $\bar{\BPi}^*=\bar{\BPi}^T=-\BPi$ and similarly  
$\left( \bar{\BPi}^{\dagger\circ} \right)^* =-{\BPi}^{\dagger\circ}  $ 
\begin{align*}
	W_1^* &=~ -\left( W_0  \left(  {\bar{\BPi}}^{\dagger\circ} \circ \big(  V_0^* A_1^*W_0     \big)   \rule{0em}{1em} \right)  \right)^*
		 ~=~  - \left( \left( {\bar{\BPi}}^{\dagger\circ}\right)^* \circ \big(   V_0^* A_1^*W_0     \big)^*   \rule{0em}{1em} \right) W_0^*\\
		 &=~   \left( {{\BPi}}^{\dagger\circ} \circ \big(  W_0^* A_1V_0     \big)   \rule{0em}{1em} \right) W_0^*
\end{align*}

\subsection{Background: Spectral Decomposition of a Matrix} 				\label{Spectral.app} 

Just as matrix partition notation is useful in expressing eigenvalue/eigenvector relations~\req{AV}, it is also useful to understand 
diagonalization and spectral decomposition of a matrix. First observe that~\req{AV} also gives the left eigenvectors of $A$ as follows
\[
	AV~=~V\Lambda ~~~~~\Leftrightarrow~~~~~  V^{-1}A ~=~ \Lambda V^{-1}~~~~~\Leftrightarrow~~~~~  A^*V^{-*} ~=~ V^{-*}\Lambda^*, 
\]
thus rows of $V^{-1}$ (columns of $V^{-*}$) are  left eigenvectors of $A$ (right eigenvectors of $A^*$).  These equations also 
give the diagonalization  
\[
			A ~=~ V\Lambda V^{-1}. 
\]
	 which  can also be  interpretated as a rank-1 decomposition of $A$ as follows. 
	 Let $W^* := V^{-1}$ (or equivalently $W := V^{-*}$), and observe that the diagonalization can be rewritten as 
	\begin{align}
		A & = 
			\obthdtall{v_1}{\cdots}{v_n}
			\thbthd{\lambda_1}{}{}{}{\ddots}{}{}{}{\lambda_n}
			\thbodwide{w^*_1}{\vdots}{w^*_n} 											\nonumber		\\
		& = 
		\lambda_1 \thbotall{v_1} \obthwide{w^*_1} + \cdots + 
		\lambda_n  \thbotall{v_n} \obthwide{w^*_n}
		~=~ \sum_{i=1}^n \lambda_i ~v_i w^*_i .
	  \label{rank1A.eq}
	\end{align}
	 This is a rank-1 (aka dyadic) decomposition of $A$, namely into $n$ rank-1 matrices  made up of  outer products
	 of the respective columns of $V$ and $W$ scaled by the respective eigenvalue. Note that since $W^*V=I$, there are 
	 simple relationships between the columns of $V$ and $W$. 
%
%
	 \begin{itemize} 
	 	\item {\em The sets $\{v_i\}$ and $\{w_i\}$ form a reciprocal basis:} 
			Reciprocal bases\footnote{Another common term is {\em dual bases}.}  
			have the property that
			$
				v_i^* w_j ~=~ \delta_{i-j}. 
			$
			This is easily seen to be true from the following partitioning of $W^*V=I$
                            \[
                            	V^{-1}V ~=~   W^*V ~=~ 
                            		\thbodwide{w^*_1}{\vdots}{w^*_n}
                            		\obthdtall{v_1}{\cdots}{v_n} 
                            		~=~ 
                            		\thbthd{w^*_1v_1}{\cdots}{w^*_1v_n}{\vdots}{}{\vdots}{w^*_nv_1}{\cdots}{w^*_nv_n} 
				~=~ I. 
                            \]
                            The reciprocal basis is useful since it allows for writing any vector $x$ in terms of a basis $\{v_i\}$
                            by observing
                            \be
                            	x ~=~ VW^* x 
				~= \obthdtall{v_1}{\cdots}{v_n}  \thbodwide{w^*_1}{\vdots}{w^*_n}   \thbotall{x}
				=~ \sum_{i=1}^n v_i \inprod{w_i}{x} . 
			      \label{recipbasis.eq}	
                            \ee			
                            Thus the coefficients of expansion of a vector $x$ in a basis $\{v_i\}$ are the inner products 
                            $ \inprod{w_i}{x}$ of the vector with the respective elements of the reciprocal basis $\{w_i\}$. 
	 
	 \item
            	 We can obtain yet another interpretation of the action of a diagonalizable  matrix on a vector 
            	 by acting with~\req{rank1A} on any vector $x$ 
            	 \be
            	 	Ax ~=~  \sum_{i=1}^n \lambda_i ~ v_i w_i^*x ~=~  \sum_{i=1}^n \lambda_i ~ v_i \inprod{w_i}{x}  
            				~=~  \sum_{i=1}^n \lambda_i ~ P_i x
            	  \label{Axdecomp.eq}
            	\ee
            	Each matrix $P_i:=v_iw_i^*$ is a (not necessarily orthogonal)  projection operator (note: $P^2=P$) onto an individual eignsubspace of $A$.  Note that the 
            	inner product is with the corresponding {\em left} eigenvector $w_i$. This expression 
            	above is known as the spectral decomposition of a linear operator. 
	
	\item 
			Let $f$ be any function analytic in a neighborhood of the set $\{\lambda_1,\ldots,\lambda_n\}$, then 
			\[
				f(A) ~=~  \sum_{i=1}^n f(\lambda_i) ~ P_i
			\]
	
	 \end{itemize}

Note that when $A$ is Hermitian, columns of $V$ are orthonormal, $W=V$, and the projections above are orthogonal.

\subsection{Repeated Eigenvalues Case for Non-normal $A_0$} 


Now assume $A_0$ has a repeated eigenvalue $\lamb$ (multiplicity $m$) and define a basis so that $A_0$ is block-diagonal 
with the first $m$ basis elements  a basis for the eigensubspace of $\lamb$. With this basis,  $A_0$, 
$V_0$, $W_0$ and $\Lambda_0$ have a block partitioning as 
\[
	A_0 = \bbm  (A_{0})_{11} & 0\\ 0 & (A_{0})_{22} \ebm ,~
	V_0 = \bbm  (V_{0})_{11} & 0\\ 0 & (V_{0})_{22} \ebm ,~	
	W_0 = \bbm  (W_{0})_{11} & 0\\ 0 & (W_{0})_{22} \ebm ,~	
	\Lambda_0 = \bbm  \lamb I & 0\\ 0 & (\Lambda_{0})_{22} \ebm
\]
where $I$ is the $m\times m$ identity matrix. 
Now take eq.~\req{Lam1}
\[
	\Lambda_0 V_1 + A_1V_0 ~=~  V_0\Lambda_1 + V_1 \Lambda_0.
\]
If we partition all matrices conformably with the above partitions,  we get 
\begin{multline*}
	\bbm  \lamb I & 0\\ 0 & (\Lambda_{0})_{22} \ebm   
	\bbm (V_1)_{11} & (V_1)_{12} \\ (V_1)_{21}  & (V_1)_{22} \ebm
	+ \bbm (A_1)_{11} & (A_1)_{12} \\ (A_1)_{21}  & (A_1)_{22} \ebm
	\bbm  (V_{0})_{11} & 0\\ 0 & (V_{0})_{22} \ebm																\\
	=~ 
	\bbm  (V_{0})_{11} & 0\\ 0 & (V_{0})_{22} \ebm
	 \bbm (\Lambda_1)_{11} & 0 \\ 0  & (\Lambda_1)_{22} \ebm
	+ 
	\bbm (V_1)_{11} & (V_1)_{12} \\ (V_1)_{21}  & (V_1)_{22} \ebm
	\bbm  \lamb I & 0\\ 0 & (\Lambda_{0})_{22} \ebm .
\end{multline*} 
Taking the $11$ block equation we get 
\[
	\lamb I_n  ~(V_1)_{11}  ~+~ (A_1)_{11} (V_{0})_{11} ~=~(V_{0})_{11} (\Lambda_1)_{11} ~+~ (V_1)_{11} ~\lamb I_n  ,
\]
but $\lamb I_n$ commutes with any matrix, so we simply get 
\[
	(A_1)_{11} (V_{0})_{11} ~=~(V_{0})_{11} (\Lambda_1)_{11},
\]
from which we can obtain $(\Lambda_1)_{11}$ by left multiplication by $(W^*_{0})_{11}$
\[
	(\Lambda_1)_{11} ~=~ 
	(W_{0})_{11}^* (A_1)_{11} (V_{0})_{11}  ~=~ (W_{0}^* A_1 V_{0})_{11} , 
\]
where the last equality follows from the block-diagonal structure of $V_0$ and $W_0$. 
Now let's examine the meaning of this. The perturbation expansion~\req{pertexp} assumes that all the $\Lambda_i$'s are
diagonal matrices (otherwise~\req{pertexp} is not an eigenvalue/eigenvector relation). This equation states that for $(\Lambda_1)_{11}$ to be diagonal as required, we must choose $(V_0)_{11}$ so that its columns are the eigenvectors 
of $(A_1)_{11}$. Consequently, $(\Lambda_1)_{11}$ will simply be the diagonal matrix of eigenvalues of $(A_1)_{11}$. 

The meaning of the above is that only this particular choice of $V_0$ will yield the expansion~\req{pertexp}. Unlike the
non-repeated case, we also have to discover $V_0$, not just $V_1$ and $\Lambda_1$.

Now to calculate $V_1$ we can apply the psuedo-inverse formula~\req{cLpinv} (which produced~\req{SylSol} earlier), but now we pay special attention 
to the structure of $\BPi^{\dagger\circ}$. In this case, it is 
\[
	\left( \BPi^{\dagger\circ} \right)_{ij} ~:= 
		\left\{ 	\begin{array}{ll}		\displaystyle	\frac{1}{\lambda_{0i} - \lambda_{0j} } , ~~~~ &\lambda_{0i} \neq \lambda_{0j}, 			\\ 
														0															, 			&\lambda_{0i} = \lambda_{0j}. 
					\end{array}			\right.
\]
It is instructive to look at the structure of $\BPi^{\dagger\circ}$, it is of the following form 
\[
	\BPi^{\dagger\circ} ~=~ \bbm \scalebox{2}{0} &    &    \\ 
																	&  0 &   \\ 
																	&     & \ddots    				\ebm ,
\]
where the large $0$ block is $m\times m$, the diagonal is all zeros, and the remaining terms are $1/(\lambda_{0i} - \lambda_{0j})$. 
The expression for $V_1$ is now the same as~\req{SylSol}, but now we also make sure to use the $V_0$ that made $\Lambda_1$ diagonal 
above (recall that the derivation of~\req{SylSol} assumed $\Lambda_1$ to be diagonal) 
\[
	V_1 ~=~ -V_0  \left(  \BPi^{\dagger\circ} \circ \big(   W_0^*A_1V_0		 \big)   \rule{0em}{1em} \right) . 
\]
Again, to compare with existing expressions, find $v_{1k}$ using $V_1e_k$ 
\be
		v_{1k} ~=~- V_1e_k ~=~ -\sum_{\lambda_{0i}\neq \lambda_{0k}}   \frac{w_{0i}^*A_1v_{0j}}{\lambda_{0i}-\lambda_{0j}}  ~v_{0i} e_j^* ~e_k
					~=~  \sum_{\lambda_{0i}\neq \lambda_{0k}}   \frac{w_{0i}^*A_1v_{0k}}{\lambda_{0k}-\lambda_{0i}}  ~v_{0i}, 
\ee

\newpage

\bibliographystyle{IEEEtran}
\bibliography{MP.bib}

\end{document}